\documentclass[twoside,leqno,10pt]{article}

\usepackage{graphics,ifthen}
\usepackage{latexsym}
\begin{document} 
\input{latex.sty}
\input{references.sty}
\input epsf.sty

\def\ind{\stackrel{\mathrm{ind}}{\sim}}
\def\iid{\stackrel{\mathrm{iid}}{\sim}}
\def\Prodi{\mathop{{\lower9pt\hbox{\epsfxsize=15pt\epsfbox{pi.ps}}}}}
\def\prodi{\mathop{{\lower3pt\hbox{\epsfxsize=7pt\epsfbox{pi.ps}}}}}

\def\Definition{\stepcounter{definitionN}\
    \Demo{Definition\hskip\smallindent\thedefinitionN}}
\def\EndDefinition{\EndDemo}
\def\Example#1{\Demo{Example [{\rm #1}]}}
\def\EndExample{\qed\EndDemo}
\def\Category#1{\centerline{\Heading #1}\rm}
\
\def\e{\text{\hskip1.5pt e}}
\newcommand{\eps}{\epsilon}
\newcommand{\proof}{\noindent {\bf Proof:\ }}
\newcommand{\remarks}{\noindent {\bf Remarks:\ }}
\newcommand{\note}{\noindent {\bf Note:\ }}
\newcommand{\examp}{\noindent {\bf Example:\ }}
\newcommand{\Lower}[2]{\smash{\lower #1 \hbox{#2}}}
\newcommand{\ben}{\begin{enumerate}}
\newcommand{\een}{\end{enumerate}}
\newcommand{\bi}{\begin{itemize}}
\newcommand{\ei}{\end{itemize}}
\newcommand{\hp}{\hspace{.2in}}

\newtheorem{lw}{Proposition 3.1, Lo and Weng (1989)}
\newtheorem{thm}{Theorem}[section]
\newtheorem{defin}{Definition}[section]
\newtheorem{prop}{Proposition}[section]
\newtheorem{lem}{Lemma}[section]
\newtheorem{cor}{Corollary}[section]
\newcommand{\rb}[1]{\raisebox{1.5ex}[0pt]{#1}}
\newcommand{\mc}{\multicolumn}

\def\Beta{\text{Beta}}
\def\Dir{\text{Dirichlet}}
\def\DP{\text{DP}}
\def\P{{\bf p}}
\def\fhat{\widehat{f}}
\def\GA{\text{gamma}}
\def\ind{\stackrel{\mathrm{ind}}{\sim}}
\def\iid{\stackrel{\mathrm{iid}}{\sim}}
\def\K{{\bf K}}
\def\min{\text{min}}
\def\N{\text{N}}
\def\p{{\bf p}}
\def\U{{\bf U}}
\def\u{{\bf u}}
\def\w{{\bf w}}
\def\X{{\bf X}}
\def\Y{{\bf Y}}

\newcommand{\reals}{{\rm I\!R}}
\newcommand{\PR}{{\rm I\!P}}
\def\Z{{\bf Z}}
\def\yy{{\mathcal Y}}
\def\rr{{\mathcal R}}
\def\BP{\text{beta}}
\def\ts{\tilde{t}}
\def\js{\tilde{J}}
\def\gs{\tilde{g}}
\def\fs{\tilde{f}}
\def\ys{\tilde{Y}}
\def\ps{\tilde{\mathcal {P}}}

\def\Report{Lancelot F. James}
\def\Author{Markov Krein}
\pagestyle{myheadings}
\markboth{\Author}{\Report}
\thispagestyle{empty}

\bct\Heading
Functionals of Dirichlet processes, the Markov Krein Identity and
Beta-Gamma processes
\lbk\lbk\smc
Lancelot F. James\footnote{
\eightit AMS 2000 subject classifications. 
               \rm Primary 62G05; secondary 62F15.\\
\eightit Corresponding authors address.
                \rm The Hong Kong University of Science and Technology,
Department of Information Systems and Management,
Clear Water Bay, Kowloon,
Hong Kong.
\rm lancelot\at ust.hk\\
\indent\eightit Keywords and phrases. 
                \rm  
          Beta-gamma processes,
          Dirichlet process,
          Markov-Krein identity,   
          gamma process.
          }
\lbk\lbk
\BigSlant 
The Hong Kong University of Science and Technology\rm 
\ect 
\Quote
This paper describes how one can use the well-known Bayesian prior to posterior
analysis of the Dirichlet process, and less known results for the
gamma process, to address the formidable problem of assessing the
distribution of linear functionals of Dirichlet processes. In
particular, in conjunction with a gamma identity, we show easily that
a generalized Cauchy-Stieltjes transform of a linear functional of a
Dirichlet process is equivalent to the Laplace functional of a class
of, what we define as, beta-gamma processes. This represents a generalization of the
Markov-Krein identity for mean functionals of Dirichlet
processes. A prior to posterior analysis of beta-gamma processes is
given that not only leads to an easy derivation of the Markov-Krein
identity, but additionally yields new distributional identities
for gamma and beta-gamma processes. These results give new
explanations and intepretations of exisiting results in the
literature. This is punctuated by establishing a simple distributional
relationship between beta-gamma and Dirichlet processes. 
\EndQuote
\rm
\newpage

\section{Introduction}
Let $P$ denote a Dirichlet random probability measure on a Polish space
${\yy}$, with law denoted
as ${\mathcal D}(dP|\theta H)$, where $\theta$ is a non-negative scalar and
$H$ is a (fixed) probability measure on ${\yy}$. In addition
let $\mm$ denote the space of boundedly finite measures on $\yy$. This
space contains the space of probability measures on $\yy$. 
The Dirichlet
process was first made popular in Bayesian nonparametrics by Ferguson
(1973), see also Freedman (1963) for an early treatment, and has subsequently been used in numerous statistical
applications. Additionally, the Dirichlet process arises in a variety
of interesting contexts outside of statistics. Formally, $P$ is said to
be a Dirichlet process if and only if for each finite collection of disjoint
measureable sets $A_{1},\ldots, A_{k}$, the random vector $P(A_{1}),\ldots,
P(A_{k})$ has a Dirichlet distribution with parameters
$\theta H(A_{1}),\ldots, \theta H(A_{k})$. In particular $P(A)$ is a
beta random variable for any measurable set $A$. An important
representation of the Dirichlet, which is analogous to Lukacs
characterization of the gamma distribution, is 
\Eq
P(\cdot)=\frac{\mu(\cdot)}{T}
\label{gamma}
\EndEq
where $\mu$ is a gamma process with finite shape parameter $\theta H$
and $T=\int_{\yy}\mu(dy)$ is a gamma random variable with shape
$\theta$ and scale $1$. 
The law of the gamma process is denoted as ${\mathcal G}(d\mu|\theta H)$,
and is characterized by its Laplace functional
\Eq
\int_{\mm}{\mbox e}^{-\mu(g)}{\mathcal {G}}(d\mu|\theta H)={\mbox e}^{-\int_{\yy}\log[1+g(y)]\theta H(dy)}
\label{laplace}
\EndEq
for $g$ a positive measurable function satisfying,
\Eq
\int_{\yy}\log[1+g(y)]\theta H(dy)<\infty.
\label{gcond}
\EndEq
We note that if $H$ is non-atomic
then the gamma process may be represented as 
\Eq
\mu(dy)=\int_{0}^{\infty}sN(ds,dy)
\label{gampoi}
\EndEq
where $N$ is a Poisson random measure with mean intensity,
$$
\theta s^{-1}{\mbox e}^{-s}dsH(dy).
\label{gamlevy}
$$
An important fact is that
$T$ and $P$ are independent, which as we shall see, has a variety of implications. 

An interesting and formidable problem initiated in a series of papers
by Cifarelli and Regazzini(1990) is the study of the exact
distribution of linear functionals, 
$$
P(g)=\int_{\yy}g(y)P(dy)
$$
of the Dirichlet process. Diaconis and Kemperman~(1996)
discuss an important by-product of this work called the {\it
 Markov-Krein} identity for means of Dirichlet processes, 
\Eq
\int_{\mm}\frac{1}{{(1+zP(g))}^{\theta}}{\mathcal {D}}(dP|\theta H)
={\mbox e}^{-\int_{\mathcal Y}\log[1+zg(y)]\theta H(dy)},
\label{MK1}
\EndEq
which has implications relative to the Markov moment
problem, continued fractions theory, exponential representations of
analytic functions, etc [see Kerov~(1998) and Tsilevich, Vershik and Yor~(2001a)].  Tsilevich, Vershik and Yor~(2001a) expand upon this emphasizing that
the right hand side of~\mref{MK1} is the Laplace functional of a Gamma
process with shape $\theta H$. That is, 
\Eq
 \int_{\mm}\frac{1}{{(1+zP(g))}^{\theta}}{\mathcal {D}}(dP|\theta H)
=
\int_{\mm}{\mbox e}^{-z\mu(g)}{\mathcal {G}}(d\mu|\theta H)
\label{MK2}
\EndEq
Their interpretation of the Markov-Krein identity  
is that the generalized Cauchy-Stieltjes transform, of order $\theta$,
of $P(g)$, where P is a Dirichlet process with shape $\theta H$,  is the Laplace transform of $\mu(g)$
when $\mu$ is the gamma process with shape $\theta H$. The authors
then exploit this fact to rederive~\mref{MK2} via an elementary proof
using the independence property of $P$ and $T$. An interesting
question, is what can one say about
\Eq
\int_{\mm}\frac{1}{{(1+zP(g))}^{q}}{\mathcal {D}}(dP|\theta H)
\label{genMK}
\EndEq
when $\theta$ and $q$ are arbitrary positive numbers? That is, can one
establish a relationship of~\mref{genMK} to the Laplace functional of
some random measure say $\mu^{*}$, which is similar to $\mu$,  for all $q$ and $\theta$. Lijoi and
Regazzini~(2003) establish analytic results for~\mref{genMK}, relating
them to the Lauricella theory of multiple hypergeometric
functions. Theorem 5.2 of their work gives analogues
of~\mref{MK1}, stating what they call a {\it Lauricella identity},
but does not specifically state a relationship such as~\mref{MK2}. We
should say for the case $\theta>q$ that it would not be terribly
difficult to deduce an analogue of~\mref{genMK} from
their result. However, this is not the case when $\theta<q$, which is
expressed in terms of contour integrals. Their representations, for all
$\theta$ and $q$, as
clearly demonstrated by the authors, however have practical utility in
regards to formulae for the density of~$P(g)$. In this case, one wants
to have an expression for~\mref{genMK}, when $q=1$ and for all $\theta$. 
Related works include
the papers of Regazinni, Guglielmi and di Nunno~(2002), Regazzini, Lijoi
and Pruenster~(2003) and the manuscript of James~(2002).

\subsection{Preliminaries and outline}
In this paper we develop results that are complementary to
the work of Lijoi and Regazzini~(2003) and Tsilevich, Vershik and
Yor~(2001a). In particular, we show that the Markov-Krein identity, as interpreted in
Tsilevich, Verhsik and Yor~(2001a), extends to a relationship
between~\mref{genMK} and the Laplace functional of a class of what we
call beta-gamma processes defined by scaling the gamma process law by
$T^{-(\theta-q)}$, for all positive $\theta$ and $q$. 
That is processes with laws equal to 
\Eq
{\mathcal {BG}}(d\mu|\theta H,
\theta-q)=\frac{\Gamma(\theta)}{\Gamma(q)}T^{-(\theta-q)}
{\mathcal G}(d\mu|\theta H)
\label{betagamma}
\EndEq

Perhaps more interesting,  is the method of approach, and derivation of supporting
results, used to establishing
such a result, which is quite different than the analytic techniques
used previously. The approach relies in part on, in this case mostly
familar, Bayesian prior posterior
calculus for Dirichlet and gamma processes in conjunction with the
usage of the following well-known {\it gamma identity} for $q>0$
\Eq
T^{-q}=\frac{1}{\Gamma(q)}\int_{0}^{\infty}v^{q-1}{\mbox e}^{-vT}dv.
\label{gammaid}
\EndEq
That is to say, purely analtyic arguments are replaced by probabilistic
ones using the familiar results in Ferguson~(1973), Lo~(1984) and
Antoniak~(1974). Thus giving the derivations a much more interpretable
Bayesian flavor. More specifically, albeit less well known, we use the
results in Lo and Weng (1989) as demonstrated for more general
proceses in James (2002, 2003). 
This bypasses the need for instance to verify certain integrability
conditions and the usage of limiting arguments. 
Moreover, somewhat conversely to Lijoi and Regazinni (2003),
we show how properties of the Dirichlet and beta-gamma processes 
yields easily interesting identities related to Lauricella and
Liouville integrals[see Lijoi and Regazinni~(2003) and Gupta and
Richards~(2001)].  
   Although we exploit the independence property of $T$ and $P$ to
prove our results, our approach is quite different from the methods
used  by Tsilevich, Verhsik and Yor (2001) to
prove~\mref{MK2}. While their proof is certainly elegant, it does not
seem possible to extend to other processes. Our methods however are
influenced by their proof of an analogous result for the two-parameter
extension of the Dirichlet process which relies on~\mref{gammaid} and
the fact that such processes are based on scaled laws. That is to say
we present an approach which is extendable to other models[see James
(2002, section 6)]. 

The study of scaled laws are of clear interest in
the case of the stable law of index $0 <\alpha<1$ as discussed in
Pitman and Yor (1997)[see also Pitman (2002)]. In particular the
two-parameter $(\alpha, \theta)$ extension of the Dirichlet process[see Pitman (1996)]
can be  defined as $P$ in~\mref{gamma} which has law governed by 
$T^{-\theta}{\mathcal P}(d\mu|\rho_{\alpha},H)$, where 
${\mathcal P}(d\mu|\rho_{\alpha},H)$ is the law of the stable
law process, which can be derived from a Poisson random measure with
intensity, 
$$
\rho_{\alpha}(ds)H(dx)=s^{-\alpha-1}dsH(dx).
$$
However, for the beta-gamma processes defined in~\mref{betagamma},
the independence property between $T$ and $P$ translates into the
following property
\Eq
\int_{\mm}f(P){\mathcal D}(dP|\theta H)=
\int_{\mm}f(P){\mathcal G}(d\mu|\theta H)=
\int_{\mm}f(P){\mathcal {BG}}(d\mu|\theta H,\theta-q).
\label{subtle}
\EndEq
for all integrable $f$.
The property~\mref{subtle} seems to suggest that the beta-gamma process may not
have much utility relative to calculations involving $P$, however it
is precisely this property that we shall exploit. 
The outline of this
paper is as follows; in section 2 we revisit well-known properties of
the Dirichlet process $P$ and describe the posterior
distribution of $\mu$ given $\Y=\{Y_{1},\ldots, Y_{n}\}$ when $\Y|P$ is
an iid sample from $P$. In particular, we show that the posterior
distribution of $\mu|\Y$ is a posterior beta-gamma process, and derive
its Laplace functional. In section 3 we present a general analysis of
the properties of the beta-gamma process including a description of
its posterior distribution which is shown to be a conjugate
family. Additionally formulae for the Laplace functional are derived
and some key
results which are relevant to the proof of the Markov-Krein identity
are given.
In sections 4 and 5 we present results that describes formally the
various relationships between~\mref{genMK} and the beta-gamma
processes.  

\section{Some results for the Dirichlet and Gamma process}
We first recall some key features of the gamma and Dirichlet
processes. Hereafter, we will assume that $H$ is non-atomic, and that
each function $g$ satisfies~\mref{gcond}.
Let $Y_{1},\ldots, Y_{n}$ denote random elements in
the space $\yy$,  which
conditional on $P$ are iid with law $P$. $P$ is a Dirichlet process
with shape $\theta H$. From Lo and Weng (1989)[see
also James (2002, 2003)], one has the following disintegration of
measures, 
\Eq
\prod_{i=1}^{n}\mu(dY_{i}){\mathcal G}(d\mu|\theta H)
={\mathcal G}(d\mu|\theta H +\sum_{i=1}^{n}\delta_{Y_{i}})
\theta^{n(\p)}\prod_{j=1}^{n(\p)}\Gamma(e_{j,n})H(dY^{*}_{j}),
\label{disint1}
\EndEq
where ${\mathcal G}(d\mu|\theta H +\sum_{i=1}^{n}\delta_{Y_{i}})$
denotes a gamma process with shape $\theta
H+\sum_{i=1}^{n}\delta_{Y_{i}}$. The quanitity $\p=\{C_{1},\ldots,
    C_{n(\p)}\}$ denotes a partition of the integers $\{1,\ldots, n\}$,
      with $n(\p)$ elements. The $C_{j}=\{i:Y_{i}=Y^{*}_{j}\}$ for $j=1,\ldots, n(\p)$
denote the collection of values equal to each unique $Y^{*}_{j}$, for
$j=1,\ldots, n(\p)$. The size of each cell $C_{j}$ is denoted as
$e_{j,n}$. Recall that a gamma process can also be described as
in~\mref{gampoi} in terms of a Poisson random measure with mean
intensity~\mref{gamlevy}. 
Another important property of the Gamma process is the
following exponential change of measure formula, 
\Eq
{\mbox e}^{-\mu(g)}{\mathcal G}(d\mu|\theta H)
={\mathcal G}(d\mu|\rho_{1+g},\theta H)
{\mbox e}^{-\int_{\yy}\log[1+g(y)]\theta H(dy)}
\label{xchange}
\EndEq
where ${\mathcal G}(d\mu|\rho_{1+g},\theta H)$ denotes a weighted gamma
process defined for each Borel measurable set $A$ as, 
\Eq
\int_{A}[1+g(y)]^{-1}\mu(dy)
\label{wg}
\EndEq
or equivalently by its inhomogeneous L\'evy measure 
$
\rho_{1+g}(ds|y)={\mbox e}^{-(1+g(y))s}s^{-1}ds. 
$ 
\Remark
The result~\mref{xchange} appears in Lo and Weng (1989, Proposition
3.1) and has been independently derived in Tsilevich, Vershik and
Yor (2001b). Versions of this result also hold for the case where $g$ is
negative or complex valued. For this study, we can actually bypass the
explicit usage of~\mref{xchange} but use it to demonstrate how one
might apply a similar result for more general random processes as in
James (2002). In
particular, we note that this operation is useful for the stable law,
which does not admit a result similar to~\mref{disint1}
\EndRemark

As an illustration, similar to Lo and Weng (1989, Corollary 3.1),  we show how one
uses~\mref{disint1} to obtain the classic results for the Dirichlet process via its
representation as a normalized gamma process. 
In particular, we 
show how one easily estabishes the disintegration, 
\Eq
\prod_{i=1}^{n}P(dY_{i}){\mathcal D}(dP|\theta H)=
{\mathcal D}(dP|\theta H+\sum_{j=1}^{n}\delta_{Y_{i}})
\pi(\p|\theta)
\prod_{j=1}^{n(\p)}H(dY^{*}_{j}),
\label{Ddisint}
\EndEq
where ${\mathcal D}(dP|\theta H+\sum_{j=1}^{n}\delta_{Y_{i}})$ is a
Dirichlet process with shape $\theta H+\sum_{j=1}^{n}\delta_{Y_{i}}$, and
is the posterior distribution of $P|\Y$ [see Ferguson (1973)]. The
partition probability
$$
\pi(\p|\theta)=\frac{\Gamma(\theta)}{\Gamma(\theta+n)}
\theta^{n(\p)}\prod_{j=1}^{n(\p)}\Gamma(e_{j,n}),
$$
is a variant of Ewens sampling formula derived by Ewens (1972) and
Antoniak (1974). The marginal distribution of $\{Y_{1},\ldots, Y_{n}\}$
is the Blackwell-MacQueen P\'olya urn distribution which can be
represented as
\Eq
\PR(d\Y|\theta H)=\pi(\p|\theta)\prod_{j=1}^{n(\p)}H(dY^{*}_{j}).
\label{Black}
\EndEq
The result in~\mref{Ddisint} can be obtained by working with the joint probability
measure of $\Y$ and $\mu$ given by replacing ${\mathcal
  D}(dP|\theta H)$
on the left hand side of~\mref{Ddisint} with the gamma process law, ${\mathcal
  G}(d\mu|\theta H)$. The task then shifts to finding the disintegration of the law of $\{\Y,\mu\}$ in terms of the posterior
distribution of $\mu|\Y$ and the marginal distribition of $\Y$. An
application of~\mref{disint1} shows that the joint distribution of $\{\Y,\mu\}$ has
the following disintegrations,
\Eq
T^{-n}\prod_{i=1}^{n}\mu(dY_{i}){\mathcal G}(d\mu|\theta H)
=T^{-n}{\mathcal G}(d\mu|\theta H +\sum_{i=1}^{n}\delta_{Y_{i}})
\theta^{n(\p)}\prod_{j=1}^{n(\p)}\Gamma(e_{j,n})H(dY^{*}_{j}).
\label{jointgamd}
\EndEq
First it is clear that
the marginal distribution of $\Y$ is obtained by integrating out
$T^{-n}$ with respect to  ${\mathcal G}(d\mu|\theta H
+\sum_{i=1}^{n}\delta_{Y_{i}})$. 
Under this law, $T$ is a gamma random variable with shape parameter
$\theta+n$ and unit scale. This implies that,
$$
\int_{\mathcal M}T^{-n}{\mathcal G}(d\mu|\theta H
+\sum_{i=1}^{n}\delta_{Y_{i}})=\frac{\Gamma(\theta)}{\Gamma(\theta+n)},
$$
yielding the desired expression for the marginal distribution of $\Y$
in~\mref{Black}.
Now in order to obtain the posterior distribution of $P$ given $\Y$,
first notice that the 
posterior distribution of $\mu$ given $\Y$ is
obtained by dividing~\mref{jointgamd} by the marginal distribution of
$\Y$. That is, we see that its posterior distribtion is of the form
in~\mref{betagamma} with specific law
\Eq
{\mathcal {BG}}(d\mu|\theta H+\sum_{i=1}^{n}\delta_{Y_{i}},n) =\frac{\Gamma(\theta+n)}
{\Gamma(\theta)}T^{-n}{\mathcal G}(d\mu|\theta H
+\sum_{i=1}^{n}\delta_{Y_{i}})
\label{post1}
\EndEq

It is not immediately obvious that the posterior distribution
of $\mu|\Y$ in~\mref{post1} indicates that the
posterior distribution $P|\Y$ is ${\mathcal D}(dP|\theta
H+\sum_{i=1}^{n}\delta_{Y_{i}})$. However, note that subject to the
gamma process law ${\mathcal G}(d\mu|\theta H
+\sum_{i=1}^{n}\delta_{Y_{i}})$, $P$ has
the desired Dirichlet process distribution and is independent of $T$. Hence for
any integrable $f$, using the independence of $P$ and $T$, it follows
from ~\mref{subtle} that the posterior distribution of $P$ given $\Y$ is characterized by
$$
\int_{\mathcal M} f(P)\pi(d\mu|\Y)=\int_{\mathcal M}f(P)
{\mathcal
  G}(d\mu|\theta H+\sum_{i=1}^{n}\delta_{Y_{i}})=\int_{\mathcal M}f(P){\mathcal
  D}(dP|\theta H+\sum_{i=1}^{n}\delta_{Y_{i}})
$$
for $f$ arbitrary integrable functions, which yields the desired result.

\subsection{The posterior distribution of the gamma process}
While the posterior distribution of $P$ given $\Y$ is well-known to be
a Dirichlet process, the corresponding posterior distribution of $\mu$
given $\Y$ described is not a Gamma process. In fact it is no longer a
L\'evy process. In order to understand what this distribution is we
find it quite useful to employ the gamma identity~\mref{gammaid}. 
In general, the gamma identity allows one to
employ an exponential change of measure formula such
as~\mref{xchange} to handle terms like $T^{-q}$. For the present setting note
that if $T$ itself is a gamma random variable with shape parameter $\tau$ and unit
scale, then,
$$
E[{\mbox e}^{-vT}]=(1+v)^{-\tau}.
$$
Now if $\tau-q>0$, Fubini's theorem yields, 
$$
E[T^{-q}]=\frac{1}{\Gamma(q)}\int_{0}^{\infty}
v^{q-1}{(1+v)}^{-\tau}dv=\frac{\Gamma(\tau-q)}{\Gamma(\tau)}
$$
This leads to the identification of a {\it gamma-gamma} density with
parameters $\tau$ and $q$ denoted as,
\Eq
{\gamma}(dv|\tau,q)=
\frac{\Gamma(\tau)}{\Gamma(q)\Gamma(\tau-q)}v^{q-1}(1+v)^{-\tau}dv
{\mbox { for }}0<v<\infty.
\label{gamgam}
\EndEq
Using the transformation $u=1/(1+v)$ yields the density of a beta random
variable with
parameters $\tau-q$ and $q$ denoted as 
\Eq
{\mathcal {B}}(du|\tau-q,q)
= \frac{\Gamma(\tau)}{\Gamma(q)\Gamma(\tau-q)}
u^{\tau-q-1}(1-u)^{q-1}
{\mbox { for }}0<u<1.
\EndEq
In this section we will encounter the special case of
$\tau=\theta+n$ and $q=n$, yielding the densities
${\gamma}(dv|\theta+n,n)$ and 
$$
{\mathcal {B}}(du|\theta,n)=
\frac{\Gamma(\theta+n)}{\Gamma(n)\Gamma(\theta)}
u^{\theta-1}(1-u)^{n-1}du
{\mbox { for }}0<u<1
$$
The gamma identity~\mref{gammaid}, ~\mref{xchange}
and Fubini's theorem yields the following explicit description of the posterior
law of $\mu|\Y$,
\begin{prop}
Let $\mu$ denote a gamma process with law ${\mathcal G}(d\mu|\theta H)$
and let $P(\cdot)=T^{-1}\mu(\cdot)$, where $T=\mu({\mathcal Y})$, denote a
Dirichlet process with shape $\theta H$. Suppose that 
$Y_{1},\ldots, Y_{n}|P$ are iid $P$, then the posterior distribution of 
$\mu|\Y$ is ${\mathcal {BG}}(d\mu|\theta H+\sum_{i=1}^{n}\delta_{Y_{i}},n)$
as defined in~\mref{post1}. Its
Laplace functional is expressible as  
\Eq
\int_{0}^{\infty}{\mbox e}^{-\int_{\mathcal Y}\log[1+\frac{1}{1+v}g(y)]\theta H(dy)}
\prod_{j=1}^{n(\p)}{(1+\frac{1}{1+v}g(Y^{*}_{j}))}^{-e_{j,n}}
{\gamma}(dv|\theta+n,n)
\label{laplaceBG}
\EndEq
The Laplace functional can also be represented in terms of a beta
density as,
\Eq
\int_{0}^{1}{\mbox e}^{-\int_{\mathcal Y}\log[1+ug(y)]\theta H(dy)}
\prod_{j=1}^{n(\p)}{(1+ug(Y^{*}_{j}))}^{-e_{j,n}}{\mathcal {B}}(du|\theta,n)
\label{laplaceBGg}
\EndEq

The Laplace functional in~\mref{laplaceBGg} shows that the posterior distribution of
$\mu|\Y$ is equivalent to the distribution of the random measure 
$U_{n}G^{*}_{n}$, where $U_{n}$ is a beta random variable with parameters 
$(\theta, n)$ and independent of $U_{n}$, $G^{*}_{n}$ is a gamma
process with shape $\theta
H+\sum_{i=1}^{n}\delta_{Y_{i}}$. As shown by~\mref{laplaceBG}, 
$U_{n}=1/(1+V_{n})$, where $V_{n}$ is a {\it gamma-gamma} random variable with
parameters $(\theta+n,n)$. 
\end{prop}
\Proof
First using~\mref{gammaid}, 
one can write
\Eq
\Gamma(n)\int_{\mm}{\mbox e}^{-\mu(g)}T^{-n}{\mathcal
  G}(d\mu|\theta H+\sum_{i=1}^{n}\delta_{Y_{i}})
=\int_{0}^{\infty}v^{n-1}\int_{\mm}{\mbox e}^{-\mu(g)}
{\mbox e}^{-vT}{\mathcal G}(d\mu|\theta H
+\sum_{i=1}^{n}\delta_{Y_{i}})dv
\label{gampf1}
\EndEq
Recalling that the exponential change of measure~\mref{xchange} transforms a gamma
process to a weighted gamma process, two applications 
show that ${\mbox e}^{-\mu(g)}
{\mbox e}^{-vT}{\mathcal G}(d\mu|\theta H
+\sum_{i=1}^{n}\delta_{Y_{i}})$ is equivalent to, 
\Eq
{\mathcal G}(d\mu|\rho_{[1+v+g]},\theta H+\sum_{i=1}^{n}\delta_{Y_{i}})
{\mbox e}^{-\int_{\mathcal Y}\log[1+\frac{1}{1+v}g(y)]\theta H(dy)}
\prod_{j=1}^{n(\p)}{(1+\frac{1}{1+v}g(Y^{*}_{j}))}^{-e_{j,n}}{(1+v)}
^{-\theta+n}
\label{lapdis}
\EndEq
Now substituting the quantities depending on $\mu$ in the right hand
side of ~\mref{gampf1} with~\mref{lapdis} yields the first expression for
the Laplace functional. The second expression is obtained by the
change of variable $u=1/(1+v)$
\EndProof 

We can now use the explicit description of the posterior distribution
of $\mu|\Y$ to deduce the posterior
distribution of $P|\Y$. Note that subject to the posterior law of
$\mu|\Y$, the posterior distribution of $P|\Y$ must be equivalent to the
law of the process
$$
P^{*}_{n}(\cdot)=\frac{U_{n}G^{*}_{n}(\cdot)}{U_{n}G^*_{n}({\mathcal Y})}
=\frac{G^{*}_{n}(\cdot)}{G^{*}_{n}({\mathcal Y})}
$$
which shows, as desired, that $P^{*}_{n}$ is a Dirichlet process with
shape $\alpha H+\sum_{i=1}^{n}\delta_{X_{i}}$.

\Remark
Although stated a bit differently, Proposition 2.1 is equivalent to
Proposition 5.7 of James (2002). 
\EndRemark
 
\section{Properties of Beta-Gamma processes}
As seen, the posterior distribution of $\mu|\Y$ is a special case of a class of
mixed gamma processes defined by a scaling operation in~\mref{betagamma}
and called beta-gamma processes. 
For shorthand, we say that $\mu$ is 
${\mathcal {BG}}(\theta H, \theta-q)$, to indicate that $\mu$
  has a beta-gamma law  defined by ${\mathcal
  {BG}}(d\mu|\theta H, \theta-q)$. When
  $\theta-q=0$, the process reduces to a gamma process with shape
  $\theta H$. Similar to the proof of Proposition 2.1, in the  case that
  $\theta-q>0$, one may use~\mref{gammaid}
  and~\mref{xchange} to obtain the Laplace functional given by the
  following proposition.
\begin{prop}
Let $\mu$ denote a beta-gamma process with law ${\mathcal
  {BG}}(d\mu|\theta H, \theta-q)$ such that $\theta-q>0$ then the
Laplace functional is given by 
\Eq
\int_{\mm}{\mbox e}^{-\mu(g)}{\mathcal {BG}}(d\mu|\theta H,
  \theta-q)
=\int_{0}^{1}{\mbox e}^{-\int_{\mathcal Y}\log[1+ug(y)]\theta H(dy)}
{\mathcal {B}}(du|q,\theta-q).
\label{laplaceBG3}
\EndEq
Using the transformation $v=u/(1-u)$ the Laplace functional can also be written in terms of a {\it
  gamma-gamma} density with parameters $(\theta-q,\theta)$.
\end{prop}

Next we describe the posterior distribution of the general class of
beta-gamma processes.

\begin{prop}
Let $\mu$ denote a beta-gamma process with law ${\mathcal {BG}}(d\mu|\theta H,\theta-q)$
and let $P(\cdot)=T^{-1}\mu(\cdot)$, where $T=\mu({\mathcal Y})$,  denote a
Dirichlet process with shape $\theta H$. Suppose that 
$Y_{1},\ldots, Y_{n}|P$ are iid $P$, then the posterior distribution of 
$\mu|\Y$ is ${\mathcal {BG}}(\theta
H+\sum_{i=1}^{n}\delta_{Y_{i}},\theta+n-q)$, defined by the law
\Eq
{\mathcal {BG}}(\theta
H+\sum_{i=1}^{n}\delta_{Y_{i}},\theta+n-q)=\frac{\Gamma(\theta+n)}
{\Gamma(q)}T^{-(\theta+n-q)}{\mathcal
    G}(d\mu|\theta H+\sum_{i=1}^{n}\delta_{Y_{i}})
\EndEq
If $\theta+n-q>0$, then the Laplace functional of the posterior distribution is 
\Eq
\int_{0}^{1}{\mbox e}^{-\int_{\mathcal Y}\log[1+ug(y)]\theta H(dy)}
\prod_{j=1}^{n(\p)}{(1+ug(Y^{*}_{j}))}^{-e_{j,n}}{\mathcal
  {B}}(du|q,\theta+n-q).
\label{laplaceBG2}
\EndEq
\end{prop}
\Proof
Note that the marginal distribution of $\Y$ is the Blackwell-MacQueen
urn distribution. The posterior distribution can be obtained easily by
following the arguments for the case of the gamma process. In
particular one can work with $T^{-(\theta+n-q)}$ in place of $T^{-n}$ in
~\mref{jointgamd} and adjust for normalizing constants. If
$\theta+n-q>0$, then the Laplace functional
is obtained by using the proof of Proposition 2.1 
with $\Gamma(n)$ and $v^{n-1}$ in~\mref{gampf1} replaced by
$\Gamma(\theta+n-q)$ and $v^{\theta+n-q-1}$.
\EndProof

Using standard Bayesian arguments, we obtain an identity for the Laplace
functional of a gamma process, and also an expression of the Laplace
functional for all beta-gamma proceses.
\begin{prop}
Let $\mu$ be a beta-gamma process, ${\mathcal {BG}}(\theta H,
\theta-q)$ then choosing an integer $n\ge 0$ such that $\theta+n-q>0$,
the Laplace functional, $\int_{\mm}{\mbox e}^{-\mu(g)}{\mathcal
  {BG}}(d\mu|\theta H, \theta -q)$, can be expressed as follows, 
$$
\sum_{\p}\pi(\p|\theta)\int_{0}^{1}{\mbox e}^{-\int_{\mathcal
    Y}\log[1+ug(y)]\theta H(dy)}\[\prod_{j=1}^{n(\p)}\int_{\mathcal
  Y}{(1+ug(y))}^{-e_{j,n}}H(dy)\]{\mathcal {B}}(du|q,\theta+n-q)
$$
for all positive $\theta$ and $q$. 
Using a change of variable $v=u/(1-u)$ leads to an expression in terms
of a {\it gamma-gamma} density. 
\end{prop}
\Proof
The proof follows from standard Bayesian
theory and the expression in~\mref{laplaceBG2}. That is, the Laplace
functional can be obtained by integrating the Laplace functional of
the posterior distribution described in~\mref{laplaceBG2} with respect
the Blackwell-MacQueen distribution as follows, 
\Eq
\int_{\yy^{n}}\[\int_{\mm}{\mbox e}^{-\mu(g)}{\mathcal {BG}}(d\mu|\theta
H+\sum_{i=1}^{n}\delta_{Y_{i}},\theta+n-q)\]\PR(d\Y|\theta H)
\label{BGlap}
\EndEq
\EndProof

We close this section with another interesting result which will be
used in the coming section.
\begin{prop}
Let $\theta,q$ be arbitrary non-negative scalars. Then for any integer
$n\ge 0$ that satisfies the constraint,
$\theta+n-q>0$ the following formula holds,  
\Eq
\int_{\mm}\frac{1}{{(T+\mu(g))}^{q}}{\mathcal G}(d\mu|\theta
H+\sum_{i=1}^{n}\delta_{Y_{i}})
=\frac{\Gamma(\theta+n-q)}{\Gamma(q)}\int_{\mm}{\mbox e}^{-\mu(g)}T^{-(\theta+n-q)}{\mathcal G}(d\mu|\theta
H+\sum_{i=1}^{n}\delta_{Y_{i}})
\label{Krein1}
\EndEq
The expressions are equivalent to 
\Eq
\frac{\Gamma(\theta+n-q)}{\Gamma(\theta+n)}\int_{\mm}{\mbox
  e}^{-\mu(g)}{\mathcal {BG}}(d\mu|\theta
H+\sum_{i=1}^{n}\delta_{Y_{i}},\theta+n-q),
\label{Kreinid}
\EndEq
whose explicit expression is deduced from~\mref{laplaceBG2}.
\end{prop}
\Proof
From~\mref{Kreinid} and~\mref{laplaceBG2}, the proof proceeds by showing that, 
\Eq
\frac{\Gamma(\theta+n)}{\Gamma(\theta+n-q)}\int_{\mm}\frac{1}{{(T+\mu(g))}^{q}}{\mathcal G}(d\mu|\theta
H+\sum_{i=1}^{n}\delta_{Y_{i}})
\label{postg}
\EndEq
is equivalent to~\mref{laplaceBG2}. Apply the gamma identity to
${(T+\mu(g))}^{-q}$ and then two applications of~\mref{xchange} to 
$
{\mbox e}^{-v[\mu(g)+T]}{\mathcal G}(d\mu|\theta
H+\sum_{i=1}^{n}\delta_{Y_{i}})
$
to show that~\mref{postg} is equal to,
$$
\frac{\Gamma(\theta+n)}{\Gamma(\theta+n-q)\Gamma(q)}
\int_{0}^{\infty}v^{q-1}{\mbox e}^{-\int_{\mathcal Y}\log[1+\frac{v}{1+v}g(y)]\theta H(dy)}
\prod_{j=1}^{n(\p)}{(1+\frac{v}{(1+v)}g(Y^{*}_{j}))}^{-e_{j,n}}{(1+v)}^{-\theta+n}dv.
$$
The result is obtained by applying the transformation $u=v/(1+v)$.
\EndProof

\section{Functionals of Dirichlet processes, The Markov-Krein
Identity and Beta-Gamma processes}
Here we show that results relating~\mref{genMK} to the beta-gamma
process can be deduced from~\mref{gammaid},~\mref{xchange} and the 
the posterior distribution of $\mu$. A general strategy
is formed by first writing 
$$
\frac{1}{1+zP(g)}=\frac{T}{T+z\mu(g)}
$$
Before proceeding to the main result we illustrate the idea for the
special case of $\theta-q>0$.
First write, 
\Eq
\int_{\mm}\frac{1}{{(1+zP(g))}^{q}}{\mathcal {D}}(dP|\theta H)
=\int_{\mm}\frac{T^{q}}{{(T+z\mu(g))}^{q}}{\mathcal {G}}(d\mu|\theta H)
\label{speccase}
\EndEq
Note the presence of $T^{q}$ causes some difficulties. However 
when $\theta-q>0$ we can use directly the special relationship between the Dirichlet
process and the beta-gamma processes exhibited
in~\mref{subtle}. That is replacing ${\mathcal {G}}(d\mu|\theta
H)$
with ${\mathcal {BG}}(d\mu|\theta H,q)$ in~\mref{speccase} yields, 
\Eq
\int_{\mm}\frac{1}{{(1+zP(g))}^{q}}{\mathcal {D}}(dP|\theta H)
=\frac{\Gamma(\theta)}{\Gamma(\theta-q)}\int_{\mm}\frac{1}{{(1+z\mu(g))}^{q}}{\mathcal
  G}(d\mu|\theta H)
\label{gaminv}
\EndEq

At this point one can evaluate the expression~\mref{gaminv}using~\mref{gammaid}. 
However a direct appeal to Proposition 3.4 with $n=0$ shows that
$$
\int_{\mm}\frac{1}{{(1+zP(g))}^{q}}{\mathcal {D}}(dP|\theta H)
=\frac{\Gamma(\theta)}{\Gamma(q)}\int_{\mm}{\mbox e}^{-\mu(g)}
T^{-(\theta-q)}{\mathcal G}(d\mu|\theta H)
=\int_{\mm}{\mbox e}^{-\mu(g)}{\mathcal {B}}(d\mu|\theta H,\theta-q)
$$

The case where $q$ is arbitrary requires another strategy that
involves~\mref{disint1} and the posterior distribution of $\mu|\Y$. 
Using these arguments
we present a new result which relates the generalized Cauchy-Stieltjes
transform of Dirichlet process linear functionals to the Laplace functional of beta-gamma processess. This presents a generalization of the
Markov-Krein identity, complementary to the Lauricella identities
deduced in Lijoi and Regazzini (2003, Theorem 5.2). 
We also present some interesting additional identities.

\begin{thm}
Let ${\mathcal D}(dP|\theta H)$ denote a Dirichlet process with shape
$\theta H$. Let $g$ denote a function satisfying~\mref{gcond}, then the
following relationships are established, 
\Enumerate
\item[(i)]
for any positive $q$ and $\theta$,
\Eq
\int_{\mm}{(1+zP(g))}^{-q}{\mathcal {D}}(dP|\theta H)=
\int_{\mm}{\mbox e}^{-z\mu(g)}{\mathcal {BG}}(d\mu|\theta H,\theta-q).
\label{MK}
\EndEq
Statement (i) implies the following results.
\item[(ii)]
For any positive $q$ and $\theta$, and any integer $n\ge 0$ which
satisfies $\theta+n-q>0$, the quantities in~\mref{MK} are equivalent to, 
$$
\sum_{\p}\pi(\p|\theta)\int_{0}^{1}{\mbox e}^{-\int_{\mathcal
    Y}\log[1+ug(y)]\theta H(dy)}\[\prod_{j=1}^{n(\p)}\int_{\mathcal
  Y}{(1+ug(y))}^{-e_{j,n}}H(dy)\]{\mathcal B}(du|q,\theta+n-q),
$$
for the gamma process, its Laplace functional may be represented as
above for all $n\ge 1$ and $q=\theta$. 
\item[(iii)]
When, $\theta-q>0$ statement (i) combined with Lemma 3.1 imply that 
\Eq
\int_{\mm}{(1+zP(g))}^{-q}{\mathcal {D}}(dP|\theta H)
=\int_{0}^{1}{\mbox e}^{-\int_{\mathcal Y}\log[1+ug(y)]\theta H(dy)}
{\mathcal {B}}(du|q,\theta-q),
\EndEq
which coincides with the result in Lijoi and Regazzini~(2003, Theorem
5, equation (5.2)). 
\EndEnumerate
\end{thm}

\Proof
For the proof of statement (i), we first assume without loss of generality
that $q=n+d$, where $d$ is a positive scalar such that
$\theta-d>0$, and $n\ge 0$ is an integer chosen such that
$\theta+n-q>0$. This means that $T^{q}=T^{n+d}$. 
Now using~\mref{subtle} with ${\mathcal {BG}}(d\mu|\theta H,d)$
yields the expression,
\Eq
\int_{\mm}{(1+zP(g))}^{-q}{\mathcal {D}}(dP|\theta H)  
=\frac{\Gamma(\theta)}{\Gamma(\theta-d)}
\int_{\mm}\frac{T^{n}}{{(T+z\mu(g))}^{q}}{\mathcal G}(d\mu|\theta H).
\label{ggproof}
\EndEq
Now to handle the term $T^{n}$, write it as
$T^{n}=\int_{\yy^{n}}\prod_{i=1}^{n}\mu(dY_{i})$, then apply~\mref{disint1}
to show that the expressions in~\mref{ggproof} are equal to, 
$$
\frac{\Gamma(\theta+n)}{\Gamma(\theta-d)}
\int_{\yy^{n}}\[\int_{\mm}\frac{1}{{(T+z\mu(g))}^{q}}{\mathcal
  G}(d\mu|\theta H+\sum_{i=1}^{n}\delta_{Y_{i}})\]
\PR(d\Y|\theta H)
$$
Apply Proposition 3.4 to the inner term, recalling that $\theta+n-q
=\theta-d$. This yields the desired expression as in~\mref{BGlap}.
\EndProof

We now discuss some interesting results obtained from Theorem
4.1. Let ${\mathcal L}(Z)$ denote the law of a random element $Z$. Recall that for a Dirichlet process functional, $P(g)$, based on
a Dirichlet process with shape $\theta H$, one can
represent its distribution as follows
\Eq
{\mathcal L}(P(g))={\mathcal L}(U_{\theta,1}P(g)+(1-U_{\theta,1})g(Y))
\label{lawDP}
\EndEq
where
$P$ is ${\mathcal D}(\theta H)$, $U_{\theta,1}$ is a beta random
variable with parameters $(\theta, 1)$ and $Y$ has distribution
H. Additionally $P$,$U$,$Y$ are independent. The right hand side
of~\mref{lawDP} just follows from the fact that for fixed $Y$, it is the
posterior distribution of $P|Y$ based on $n=1$ observation. Hence the
unconditional distribution must be ${\mathcal D}(\theta H)$. That is, 
$$
{\mathcal D}(dP|\theta H)=\int_{\yy}{\mathcal D}(dP|\theta
  H+\delta_{y})H(dy).
$$
 See
Diaconis and Freedman (1999) for an interesting usage of~\mref{lawDP}.
Here we point out some related results for $\mu$. In particular, as
pointed out in Lijoi and Regazzini~(2003), for determining the density of
functionals $P(g)$, it is useful to have the expression for the
Cauchy-Stietljes transform for the case $q=1$. Theorem 4.1 shows that 
\Eq
\int_{\mm}{(1+zP(g))}^{-1}{\mathcal {D}}(dP|\theta H)=
\int_{\mm}{\mbox e}^{-z\mu(g)}{\mathcal {BG}}(d\mu|\theta H,\theta-1).
\label{MK5}
\EndEq

Thus in principle one can use the right hand side of~\mref{MK5} to
perform an appropriate inversion to deduce the distribution of $P(g)$.
See Lijoi and Regazzini (2003, Section 6) for details in that
direction. Here, similar to the case of ~\mref{lawDP} we use the 
the explicit posterior analyis of $\mu$ to characterize beta-gamma
processes with shape $\theta H$ and $q=1$ for all $\theta$. For the
remainder of this work, let $\mu_{\theta,\theta-1}$  be ${\mathcal {BG}}(\theta
H, \theta-1)$, $U_{a,b}$ denote a beta $(a,b)$ random variable, let
$T_{p}$ denote  a gamma random variable with shape
$p$ and scale 1, and let $Y_{1}$ be a random element with
distribution $H$. Additionally, let $\mu_{\theta}$ denote a gamma
process with shape $\theta H$ and assume that the variables
$\mu_{\theta}$, $U_{a,b},T_{p},Y_{1}$ are independent. Additionally
let $P_{\theta}$ denote  a Dirichlet process with shape $\theta
H$. 

\begin{prop}
Let $\mu_{\theta,\theta-1}$ be ${\mathcal {BG}}(\theta H, \theta-1)$ then
for $g$ a real valued function such that its absolute value
satisfies~\mref{gcond}, the following distributional equalities hold;
\Enumerate
\item[(i)]
for all $\theta>0$
\Eq
{\mathcal L}(\mu_{\theta,\theta-1}(g))={\mathcal
  L}(U_{1,\theta}\mu_{\theta}(g)
+U_{1,\theta}T_{1}g(Y_{1}))
\label{lawBG}
\EndEq
\item[(ii)]
for $\theta>1$,
$$
{\mathcal L}(U_{1,\theta-1}\mu_{\theta}(g))
={\mathcal
  L}(U_{1,\theta}\mu_{\theta}(g)+
U_{1,\theta}T_{1}g(Y_{1}))
$$
\item[(iii)]
For $\theta=1$, $\mu_{1,0}:=\mu_{1}$ is a gamma process with shape $H$
and 
$$
{\mathcal L}(\mu_{1}(g))={\mathcal
  L}(U_{1,1}\mu_{1}(g)
+U_{1,1}T_{1}g(Y_{1}))
$$
\item[(iv)]
If $\mu_{\theta}$ denotes a gamma process with arbitrary shape parameter
$\theta H$ then,
\Eq
{\mathcal L}(\mu_{\theta}(g))={\mathcal
  L}(U_{\theta,1}\mu_{\theta}(g)
+
U_{\theta,1}T_{1}g(Y_{1}))
\label{lawBG4}
\EndEq
\EndEnumerate 
\end{prop}
\Proof
The proof is immediate from statement (ii) of Theorem 4.1 or
Proposition 3.3. 
One may
also use Proposition 2.1 and Proposition 3.1 to obtain statement (iv)
and (ii). It is well-known that the condition~\mref{gcond} is
necessary and sufficient for the existence of $P(g)$.
\EndProof

The equation~\mref{lawBG} provides a simple description of the law
associated with the right hand side of~\mref{MK5} and should prove
particularly useful in the case where $0<\theta<1$. Notice
that the only unknown quantity in~\mref{lawBG} is the distribution of
$\mu_{\theta}(g)$. This suggests that in principle one can concentrate
on the law of $\mu_{\theta}(g)$ to ascertain the law of $P_{\theta}(g)$ for all
$\theta$. Note this is not an obvious fact as one can use the
representation 
\Eq
T_{\theta}P_{\theta}(g)=\mu_{\theta}(g),
\label{gammaDP}
\EndEq
and the independence of $T_{\theta}$ and $P_{\theta}$,
to find easily the distribution of $\mu_{\theta}$. However if one has the
distribution of $\mu_{\theta}(g)$ one has to negotiate the typically complex dependence structure
of $T_{\theta}$ and $\mu_{\theta}$ to obtain the law of $P_{\theta}(g)$ 

One may
use ~\mref{MK5} and~\mref{lawBG} to recover the expressions in Propositions 4 and 5 of Regazzini, Guglielmi and
Di Nunno~(2002) where Proposition 5 is a special case of the results in
Regazzini, Lijoi and Pruenster (2003).
The statement, 
~\mref{lawBG4}
provides a distributional identity for an arbitrary gamma process
which we believe is new. One may verify~\mref{MK5}, by noting that~\mref{lawBG}
satisfies, 
\Eq
{\mathcal L}(\mu_{\theta,\theta-1}(g))={\mathcal
  L}(U_{1,\theta}\mu_{\theta}(g)
+
U_{1,\theta}T_{1}g(Y_{1}))={\mathcal L}(T_{1}P_{\theta}(g))
\label{lawBG6}
\EndEq
where $T_{1}=U_{1,\theta}T_{\theta+1}$, with $U_{1,\theta}$ and
$T_{\theta+1}$ independent, is independent of $P_{\theta}(g)$, where $P_{\theta}$ is 
${\mathcal D}(\theta H)$. To see this, write 
$$
U_{1,\theta}[\mu_{\theta}(g)
+T_{1}g(Y_{1})]=U_{1,\theta}T_{\theta+1}\frac{\mu_{\theta}(g)
+T_{1}g(Y_{1})}{T_{\theta+1}}
$$
and use the classic beta gamma calculus for random variables and the
usual properties of Dirichlet,
gamma processes. ~\mref{lawBG6} combined with~\mref{MK5} yield the
the obvious result,
$$
\int_{\mm}{(1+zP(g))}^{-1}{\mathcal {D}}(dP|\theta H)=
\int_{0}^{\infty}\int_{\mm}{\mbox e}^{-ztP(g)}{\mathcal D}(dP|\theta H)
{\mbox e}^{-t}dt
$$
One may also use~\mref{lawBG4} to prove~\mref{MK1} and~\mref{MK2}. 
\section{What is a beta-gamma process?}
The
final result which is a generalization of Proposition 4.1, is derived directly from Proposition 3.2 and
Proposition 3.3. This result explains what a beta-gamma process is and
indeed implies Theorem 4.1.
\begin{thm}
Let $\mu_{\theta,\theta-q}$ have distribution ${\mathcal {BG}}(\theta
H, \theta-q)$ and let $\mu_{\theta}$ denote a gamma process with shape
$\theta H$, then for all positive $\theta$ and $q$ and an integer
$n$ chosen such that $\theta+n-q>0$ , 
the following distributional equalities hold;
\Enumerate
\item[(i)]
For all $\theta$ and $q>0$ and an integer
$n$ chosen such that $\theta+n-q>0$,
\Eq
{\mathcal L}(\mu_{\theta,\theta-q})={\mathcal
  L}(U_{q,\theta+n-q}\mu_{\theta}
+U_{q,\theta+n-q}\sum_{j=1}^{n(\p)}G_{j,n}\delta_{Y^{*}_{j}}),
\label{lawBGlast}
\EndEq
where conditional on $\p$ the distinct variables on the right hand side are
mutually  independent such that, $U_{q,\theta+n-q}$ is
beta with parameters $(q,\theta+n-q)$, $\mu_{\theta}$ is a gamma process with shape
$\theta H$, $\{G_{j,n}\}$ are independent gamma random variables
with shape $e_{j,n}$ and scale $1$, and $Y^{*}_{j}$ for $j=1,\ldots, n(\p)$
are iid $H$. The distribution of $\p$ is $\pi(\p|\theta)$.
\item[(ii)]
For $\theta-q>0$,
${\mathcal L}(\mu_{\theta,\theta-q})={\mathcal
  L}(U_{q,\theta-q}\mu_{\theta})$

Statement (i) implies the following results. 
\item[(iii)]
For all positive $\theta$ and $q$
\Eq
{\mathcal L}(\mu_{\theta,\theta-q})={\mathcal
  L}(T_{q}P_{\theta})
\label{bgdp}
\EndEq
where $T_{q}$ is a gamma random variable with shape $q$ and scale $1$
independent of $P_{\theta}$ which is a Dirichlet process with shape
$\theta H$.
\item[(iv)]
For all positive $\theta$ and $q$,
$
{\mathcal L}(T_{\theta}\mu_{\theta,\theta-q})={\mathcal
  L}(T_{q}\mu_{\theta}),
$
where $T_{\theta}$ is gamma with shape $\theta$ and scale 1,
independent of  $\mu_{\theta,\theta-q}$. Similarly $T_{q}$ and
$\mu_{\theta}$ are independent.
\EndEnumerate 
\end{thm}
\Proof
The distributional identity in $(i)$ is a direct consequence of the
mixture representation deduced from Proposition 3.2 which is
deduced from the posterior distribution in Proposition 3.1. Note that all
quantities on the right side of~\mref{lawBGlast}, including $\p$, are
random. Statement (ii) follows from Proposition 3.1. We now show that
statement(iii) follows from statement (i). Notice that,
$
T_{\theta+n}:=\m_{\theta}(\yy)+\sum_{j=1}^{n(\p)}G_{j,n},
$
is a gamma random variable with shape $\theta+n$ independent of
$U_{q,\theta+n-q}$. Moreover using the mixture representation of the
Dirichlet process derived from its posterior distribution and
$\PR(d\Y|\theta H)$, it follows that,
$$
\frac{\mu_{\theta}+\sum_{j=1}^{n(\p)}G_{j,n}\delta_{Y^{*}_{j}}}
{T_{\theta+n}}
$$
is a Dirichlet process with shape $\theta H$, independent of
$T_{\theta+n}$ and $U_{q,\theta+n-q}$.
Hence  
the right hand side of~\mref{lawBGlast} can be written as,
$$
U_{q,\theta+n-q}T_{\theta+n}P_{\theta}
$$
The result is completed by noting that $U_{q,\theta+n-q}T_{\theta+n}$
is equal in distribution to $T_{q}$. Statement (iv) follows
immediately from (iii).
\EndProof
The expression~\mref{bgdp} tells us precisely that, for all $\theta$
and $q$, a  beta-gamma
process with parameters $\theta H$ and $\theta-q$ is 
equivalent in distribution to a Dirichlet process with shape $\theta
H$, scaled by
an independent gamma random variable with shape $q$. 
Hence, using this interpretation the first result in Theorem 4.1 is an 
immediate consequence of,
$$
E[{\mbox e}^{-z\mu_{\theta,\theta-q}(g)}]=\frac{1}{\Gamma(q)}\int_{0}^{\infty}t^{q-1}\[\int_{\mm} {\mbox
e}^{-ztP(g)}{\mathcal D}(dP|\theta H)\]{\mbox e}^{-t}dt
=\int_{\mm}{(1+zP(g))}^{-q}{\mathcal {D}}(dP|\theta H).
$$
However, it is the other representations of the beta-gamma process
that yield useful explicit expressions.

\Section{Closing Remarks}
The methodology discussed here relies on three main
ingredients which are not specifically linked to the Dirichlet or
Gamma process. That is, the use of a prior to posterior analysis of
$\mu$ which manifests itself in terms of a scaled distribution, the
use of a gamma identity and the use of an exponential change of
measure. The rest relies of course on some specific features of the
Dirichlet, gamma process calculus. James (2002, sections 5 and 6)
demonstrates that using the three general techniques mentioned above, that
is elements of what he calls a {\it Poisson process partition calculus},  that
one can obtain a generalization of a Markov-Krein type relationship
for all process
$P$ defined as in~\mref{gamma} for all random measures $\mu$ which are
completely random measures or have laws based on a scaling of such
processes. In so doing, he shows a natural duality between the posterior
distribution of $\mu$ and $P$ and the distribution of functionals of
$P$. A key point is that the law of $\mu$ can always be represented in
terms of a mixture representation derived from it posterior
distribution.
In principle, this serves to answer a question raised in Tsilevich,
Vershik and Yor~(2001a). We should say that in retrospect some of
those results in James (2002, section 5 and 6) may seem a bit cryptic.
The result in Theorem 4.1, without the explicit reference to a 
beta-gamma process are contained in Proposition 6.1 and Proposition 6.2 of
James (2002, section 6). It is of interest to clarify and refine the results in
James (2002, section 6) for more general $P$.
In particular,  we believe it would be useful to combine these ideas 
with the related results of Regazzini, Lijoi and Pruenster~(2003).

\centerline{\Heading References}
\vskip0.4in
\tenrm
\def\smc{\tensmc}
\def\sl{\tensl}
\def\bf{\tenbold}
\baselineskip0.15in

\Ref   
\by    Antoniak, C. E. 
\yr    1974
\paper Mixtures of Dirichlet processes with applications to Bayesian
nonparametric problems
\jour  \AnnStat
\vol   2
\pages 1152-1174
\EndRef

\Ref
\by    Blackwell, D. and MacQueen, J. B.
\yr    1973
\paper Ferguson distributions via P\'olya urn schemes
\jour  \AnnStat
\vol   1
\pages 353-355
\EndRef

\Ref
\by    Cifarelli, D. M. and Regazzini, E. 
\yr    1990
\paper Distribution functions of means of a Dirichlet process
\jour  \AnnStat
\vol   18
\pages 429-442
\EndRef

\Ref
\by    Diaconis, P. and Freedman, D. A. 
\yr    1999
\paper Iterated random functions
\jour  Siam Rev.
\vol   41
\pages 45-76
\EndRef

\Ref
\by    Diaconis, P. and Kemperman, J. 
\yr    1996
\paper Some new tools for Dirichlet priors. Bayesian Statistics 5 
       (J.M. Bernardo, J.O. Berger, A.P. Dawid and A.F.M. Smith eds.), 
       Oxford University Press, pp. 97-106
\EndRef 

\Ref
\by    Ewens, W. J.
\yr    1972
\paper The sampling theory of selectively neutral alleles
\jour  Theor. Popul. Biol.
\vol   3
\pages 87-112
\EndRef

\Ref
\by    Ferguson, T. S.
\yr    1973
\paper A Bayesian analysis of some nonparametric problems
\jour  \AnnStat
\vol   1
\pages 209-230
\EndRef

\Ref   
\by     Freedman, D. A.
\yr     1963
\paper  On the asymptotic behavior of
  Bayes estimates in the discrete case
\jour   Ann. Math. Statist.
\vol    34
\pages  1386-1403
\EndRef

\Ref
\yr     2001
\by     Gupta, R. D.  and Richards, D. St. P.
\paper  The history of the Dirichlet and Liouville distributions 
\jour   International Statistical Review 
\vol    69
\pages  433-446
\EndRef

\Ref
\yr    2003
\by    James,~L.~F.
\paper Poisson process partition calculus with applications to\\
       Bayesian L{\'e}vy moving averages and shot-noise processes. Manuscript
\EndRef

\Ref
\by     James, L. F. 
\yr     2002 
\paper Poisson process partition calculus with applications to
exchangeable models and Bayesian nonparametrics,
arXiv:math.PR/0205093, 2002.\\
Available at http://arxiv.org/abs/math.PR/0205093
\EndRef

\Ref
\by    Kerov, S. 
\yr    1998
\paper Interlacing measures
\jour  Amer. Math. Soc. Transl
\vol   181
\pages 35-83
\EndRef

\Ref
\by    Lijoi, A. and Regazzini, E.
\yr    2003
\paper Means of a Dirichlet process and multiple hypergeometric
functions. To appear {\it Annals of Probability}. Available as 
Quaderno di Dipartmento di Economia Politica no. 138, Pavia.
Available at http://economia.unipv.it/menu.htm 
\EndRef

\Ref
\by    Lo, A. Y.
\yr    1984
\paper On a class of Bayesian nonparametric estimates: I.  Density
       Estimates
\jour  \AnnStat
\vol   12
\pages 351-357
\EndRef 

\Ref
\by    Lo, A. Y. and Weng, C. S.
\yr    1989
\paper On a class of Bayesian nonparametric estimates: II. Hazard
       rates estimates 
\jour  Ann. Inst. Stat. Math.
\vol   41
\pages 227-245
\EndRef

\Ref
\by    Pitman, J.
\yr    2002
\paper Combinatorial stochastics processes. Lecture notes for
Saint-Flour summer school, July 2002.
Available at http://stat-www.berkeley.edu/users/pitman/bibliog.html
\EndRef
 
\Ref
\by    Pitman, J. 
\yr    1996
\paper Some developments of the Blackwell-MacQueen urn scheme. In
Statistics, Probability and Game Theory T.S. Ferguson, L.S. Shapley
and J.B. MacQueen editors, IMS Lecture Notes-Monograph series, Vol 30,
pages 245-267
\EndRef

\Ref
\by    Pitman, J. and Yor, M.
\yr    1997
\paper The two-parameter Poisson-Dirichlet distribution derived from
       a stable subordinator
\jour  \AnnProb
\vol   25
\pages 855-900
\EndRef

\Ref   
\by    Regazzini, E., Guglielmi, A. and Di Nunno, G. 
\yr    2002
\paper Theory and numerical analysis for exact distributions of
functionals of a Dirichlet process
\jour  \AnnStat
\vol   30
\pages 1376-1411
\EndRef

\Ref   
\by    Regazzini, E., Lijoi, A. and Pruenster, I.
\yr    2003
\paper Distributional results for means of normalized random measures
with independent increments
\jour  \AnnStat
\vol   31
\pages 560-585
\EndRef

\Ref   
\by    Tsilevich, N. V., Vershik, A. M, and Yor, M.
\yr    2001a
\paper Remarks on the Markov-Krein identity and quasi-invariance of
       the gamma process (in Russian). Zapiski nauchnyh seminarov POMI
       No.283, pp.21-36. English translation to appear in Journal of 
       Mathematical Sciences
\EndRef

\Ref
\by    Tsilevich, N. V., Vershik, A. M, and Yor, M.
\yr    2001b
\paper An infinite-dimensional analogue of the Lebesque measure and
       distinguished properties of the gamma process 
\jour   J. Funct. Anal 
\vol    185 
\pages   274-296
\EndRef

\vskip0.75in

\smc

\Tabular{ll}

Lancelot F. James\\
The Hong Kong University of Science and Technology\\
Department of Information Systems and Management\\
Clear Water Bay, Kowloon\\
Hong Kong\\
\rm lancelot\at ust.hk\\

\EndTabular

\end{document}